\documentclass{amsart}

\usepackage[english]{babel}

\usepackage{amssymb}
\usepackage{amsmath}
\usepackage{amsthm}
\usepackage{enumerate}
\usepackage{cite}
\usepackage{xcolor}
\usepackage[final]{listings}
\usepackage[section]{algorithm}
\usepackage{booktabs}
\usepackage{multirow}

\title[Identifiability of symmetric tensors]
{On generic identifiability of symmetric tensors of subgeneric rank}
\date{}

\newcommand{\C}{\mathbb{C}}

\newcommand{\Pj}{\mathbb{P}}

\newcommand{\Var}[1]{\mathcal{#1}}
\newcommand{\vect}[1]{\mathbf{#1}}

\newcommand{\Sec}[2]{\sigma_{#1}({#2})}
\newcommand{\Tang}[2]{\mathrm{T}_{#1} {#2}}

\newcommand{\Plane}[1]{\mathrm{#1}}
\newcommand{\rank}{\operatorname{rank}}

\newcommand{\refthm}[1]{{Theorem \ref{#1}}}
\newcommand{\refeqn}[1]{{(\ref{#1})}}

\newcommand{\refsec}[1]{{section \ref{#1}}}

\newcommand{\reftab}[1]{{Table \ref{#1}}}
\newcommand{\refprop}[1]{{Proposition \ref{#1}}}
\newcommand{\reflem}[1]{{Lemma \ref{#1}}}

\newcommand{\propref}[1]{\refprop{#1}}

\newtheorem{defn0}{Definition}[section]
\newtheorem{prop0}[defn0]{Proposition}
\newtheorem{thm0}[defn0]{Theorem}
\newtheorem{lemma0}[defn0]{Lemma}
\newtheorem{coro0}[defn0]{Corollary}

\newtheorem{rem0}[defn0]{Remark}
\newcommand{\rig}{\smash{\mathop{\longrightarrow}\limits}}


\author[L.~Chiantini]{Luca Chiantini}
\address{Luca Chiantini: Dipartimento di Ingegneria dell'Informazione e Scienze Matematiche, 
Universit\`a di Siena, Italy}
\email{luca.chiantini@unisi.it}

\author[G.~Ottaviani]{Giorgio Ottaviani}
\address{Giorgio Ottaviani: Dipartimento di Matematica e Informatica ``Ulisse Dini'', 
Universit\`a di Firenze, Italy}
\email{ottavian@math.unifi.it}
\thanks{The first and second author are members of the Italian GNSAGA-INDAM}

\author[N.~Vannieuwenhoven]{Nick Vannieuwenhoven}
\address{Nick Vannieuwenhoven: Department of Computer Science, KU Leuven, Belgium}
\email{nick.vannieuwenhoven@cs.kuleuven.be}
\thanks{The third author was supported by a Ph.~D.\ fellowship of the 
Research Foundation--Flanders (FWO)}

\begin{document}

\begin{abstract}
We prove that the general symmetric tensor in $S^d\C^{n+1}$ of rank $r$ 
is identifiable, provided that $r$ is smaller than the generic rank. 
That is, its Waring decomposition as a sum of $r$ powers of linear 
forms is unique. Only three exceptional cases arise, all of which were
known in the literature.
Our original contribution regards 
the case of cubics ($d=3$), while for $d\ge 4$ we rely on known results 
on weak defectivity by Ballico, Ciliberto, Chiantini, and Mella.
\end{abstract}

\maketitle

\section{Introduction}

We denote by $S^d\C^{n+1}$ the space of symmetric tensors on $\C^{n+1}$;
such tensors can be identified with homogeneous polynomials of 
degree $d$ in $n+1$ variables, which are also referred to as
\emph{forms}.
In this symmetric setting, the most natural tensor rank decomposition is the 
classical Waring decomposition, which expresses a symmetric tensor as a sum of 
powers of linear forms. Precisely, every form $f\in S^d\C^{n+1}$ has a minimal
expression
\begin{equation} \label{eq:Waring}
f=\sum_{i=1}^r l_i^d,
\end{equation} 
where $l_i\in\C^{n+1}$ are linear forms \cite{IK1999};
the minimal number of summands $r$ is called the symmetric rank of $f$, since 
 in the correspondence between forms and symmetric tensors,
powers of linear forms correspond to tensors of rank $1$.
A natural question concerns the number of summands required for
representing a general form in $S^d \C^{n+1}$. This problem is elementary for $d=2$, which
corresponds to the case of symmetric matrices. For $d\geq 3$, the 
question was answered by Alexander and Hirschowitz in \cite{AH1995}.
Letting
\begin{align} \label{eqn_exp_rank}
r_{d,n} = \frac{ \binom{n+d}{d} }{ n+1 },
\end{align}
they proved that the general $f\in S^d\C^{n+1}$ with $d\ge 3$ has rank 
$\lceil r_{d,n} \rceil$, which is called the \emph{generic rank},
unless the space $S^d \C^{n+1}$ is one of the so-called defective cases 
$S^4\C^{n+1}$ for $n=2, 3, 4$ and $S^3\C^{5}$,
where the generic rank is $\lceil r_{d,n}\rceil+1$. 
When the rank of a Waring decomposition is strictly smaller than $r_{d,n}$, 
we say that this decomposition is of \emph{subgeneric rank}. 
It is worth noting that, in our notation, being of \emph{subgeneric rank} is 
not always equivalent to being of \emph{rank smaller than the one of a general tensor}, 
because in the defective cases
above a general tensor has rank strictly bigger than $r_{d,n}$.

The Alexander--Hirschowitz theorem implies that the generic
tensor of subgeneric rank admits only a finite number of alternative Waring 
decompositions \cite{Land_book}.
In this paper, we shall be concerned with
proving that the\footnote{Since the generic tensor of fixed subgeneric rank $r$ is an element of the $r$-secant variety of the Veronese variety $v_d(\Pj^n)$, which is an irreducible variety \cite{Zak1993}, this terminology is warranted.} generic tensor of fixed subgeneric rank admits precisely one 
Waring decomposition, modulo permutations of the summands and scaling by $d$-roots of unity.
More precisely, the main result of this paper is the following theorem.

\begin{thm0}\label{thm_main} Let $d\ge 3$. The general tensor in $S^d\C^{n+1}$
of subgeneric rank $r < r_{d,n}$ with $r_{d,n}$ as in \refeqn{eqn_exp_rank} has a \emph{unique}
Waring decomposition, i.e., it is \emph{identifiable},
unless it is one of the following cases:
\begin{enumerate}
\item $d=6$, $n=2$, and $r=9$;
\item $d=4$, $n=3$, and $r=8$; 
\item $d=3$, $n=5$, and $r=9$.
\end{enumerate}
In all of these exceptional cases, there are exactly two Waring decompositions.
\end{thm0}

The three exceptional cases were already known in the literature. 
The first two cases are classical; see Remark 4.4 in \cite{Mella2006} and 
Remark 6.5 in \cite{C2001}. The third case was recently found by Ranestad and Voisin; 
see the proof of Lemma 4.3 in \cite{Ranestad2013}. A uniform treatment of these cases is 
presented in \refprop{prop:n235}. Our original contribution establishes that there 
are no more exceptions to identifiability for cubics.

The proof of \refthm{thm_main} is based on the study of the geometric concepts 
of \emph{weak defectivity}, developed in \cite{ChCi2001}, and 
\emph{tangentially weak defectivity}, developed in \cite{BCO2013}. 
Indeed, from this point of view, the theorem can be reformulated in 
the following way, which is a result of independent interest.

\begin{thm0}
\label{thm_main2}
Let $d\ge 3$, $r_{d,n}$ as in \refeqn{eqn_exp_rank}, and $r < r_{d,n}$.
Then, the common singular locus of the space of hypersurfaces of degree 
$d$ in $\Pj^n$ that are singular at $r$ general points,
consists of exactly these $r$ points, except in the following cases:
\begin{enumerate}
\item $d=6$, $n=2$, and $r=9$. The unique sextic plane curve singular 
at $9$ general points is a double cubic, so that its singular locus is an
elliptic cubic curve;
\item $d=4$, $n=3$, and $r=8$. The net of quartic surfaces singular at 
$8$ points consists of reducible quadrics, so that the common singular locus 
is the base locus of the pencil of quadrics through $8$ general points, 
which is an elliptic normal curve of degree $4$;  
\item $d=3$, $n=5$, and $r=9$. The common singular locus of the pencil of 
cubic 4-folds singular at $9$ general points is the unique elliptic normal 
curve of degree $6$ through these $9$ points.  
\end{enumerate}
Furthermore, the above exceptional cases are the only instances where there
exists a unique elliptic normal curve of degree $n+1$ in $\Pj^n$ through
$r$ general points.
\end{thm0}

In this formulation, the theorem was already partially proved:
the case $n \le 2$ was proved by Chiantini and Ciliberto \cite{ChCi2006}; 
for $d\ge 4$ it was proved by Ballico \cite[Theorem 1.1]{BL}; 
and for $d = 3$ with $r < r_{d,n} - \frac{n+2}{3} + 1$ it was proved by Mella
\cite[Theorem 4.1]{Mella2006}. 
Consequently, the original contribution of this paper concerns the case of 
cubics, i.e., $d=3$, which we solve completely in the subgeneric case. 
This answers the question posed in Remark 4.4 in \cite{Mella2006}.

We notice that Ballico \cite{BL} proved an even stronger result for $d\ge 4$. 
Namely, he showed that a general hypersurface of degree $d\ge 4$ in $\Pj^n$ that
is singular in $r$ general points, is singular only at these $r$ points
(except for the exceptional cases (1) and (2) of \refthm{thm_main2}).
This is equivalent to showing that the Veronese variety $v_d(\Pj^n)$ is not 
$r$-weakly defective, while our result only says that it is not 
\emph{$r$-tangentially} weakly defective. We wonder whether the above list also 
gives the classification of all $r$-weakly defective Veronese varieties
$v_d(\Pj^n)$, 
even for $d=3$.

Symmetric tensors of general rank are not expected to admit only a finite number of 
Waring decompositions, because the expected dimension $\lceil r_{d,n} \rceil (n+1)$ 
of the $\lceil r_{d,n} \rceil$-secant variety of the Veronese variety $v_d(\Pj^n)$ 
may exceed the dimension $\binom{n+d}{d}$ of the embedding space $S^d \C^{n+1}$. Therefore,
at least a curve's worth of alternative Waring decompositions of a general symmetric 
tensor is anticipated in these cases. However, if $r_{d,n} = \lceil r_{d,n} \rceil$ is integer, then 
a general symmetric tensor is still expected to admit only a finite number of Waring 
decompositions. The approach pursued in this paper, i.e., proving not tangential weak 
defectivity, cannot handle tensors of the generic rank. Other approaches need to be 
considered in this setting. In fact, Mella \cite{Mella2009} formulated a conjecture 
about the cases where the expression in \refeqn{eq:Waring} is still expected to be 
unique even for general symmetric tensors. In \cite{Bertini2014}, further evidence for 
this conjecture was given; in addition, the analogous problem for nonsymmetric tensors 
was also considered.

Even though the general symmetric tensor is not of subgeneric rank, the setting 
considered in this paper is nevertheless important in applications where
one is mostly interested in the identifiability of symmetric tensors of subgeneric rank.
For instance, Anandkumar, Ge, Hsu, Kakade, and Telgarsky \cite{Anandkumar2014}
consider statistical parameter inference algorithms based on decomposing symmetric 
tensors for a wide class of latent variable models. The identifiability of the Waring
decomposition then ensures that the recovered parameters, which correspond with the 
individual symmetric rank-$1$ terms in Waring's decomposition, are unique, and, thus, 
admit an interpretation in the application domain. The rank of the Waring decomposition, 
in these applications, is invariably much smaller than the generic rank.
As general sources on tensor decomposition, we refer to \cite{Strassen1983,IK1999,C2001,Comon2008,Land_book}.

In analogy to \refthm{thm_main}, we mention that
the results in \cite{BCO2013,CO,COV2014} give broad evidence to the
analogous problem in the setting of nonsymmetric tensors, i.e., that a general 
nonsymmetric tensor of subgeneric rank admits a unique tensor rank
decomposition, unless it is one of the exceptional cases that have already 
been proved in \cite{AOP2009,CO,BCO2013,CMO14,BC13}.

In the proofs by induction of several theorems, we rely on the software Macaulay2 \cite{M2} for proving the base cases.
The two scripts we used are available as ancillary files in the arXiv submission of this paper.

The content of the paper is the following.
In \refsec{sec:newexample}, we present a uniform treatment of the exceptional 
cases appearing in Theorems \ref{thm_main} and \ref{thm_main2}. Remark \ref{rem:genesis} 
also discusses our initial motivation for studying the topic of this paper. 
Section \ref{sec:cubics} contains the proof of the main theorem.
Thereafter, the connection between weakly defective varieties and the dual
varieties to secant varieties, including a description of the dual varieties of
all weakly defective examples appearing in
\refthm{thm_main2}, is explored in \refsec{sec:dual}.
In particular, \refthm{thm:dualveronese} contains the description of 
cubic hypersurfaces in $\Pj^5$ that can be written as the determinant of a 
$3\times 3$ matrix with linear entries. In \refsec{sec:algorithm},
we give an explicit criterion allowing to check if a given Waring decomposition 
is unique. This algorithm is an extension to the symmetric case of the one
provided in \cite{COV2014} for general tensors.

\smallskip
\paragraph{{\bf Acknowledgements}}
We want to thank C.~Ciliberto, I.~Domanov, J.M.~Landsberg, M.~Mella, L.~Oeding, 
K.~Ranestad and F.~Russo for useful discussions. In particular, 
K.~Ranestad pointed out the use of Gale transforms to prove \refprop{prop:n235} 
and informed us about Lemma 4.3 of \cite{Ranestad2013} which contains a proof 
of the third exceptional case in Theorems \ref{thm_main} and \ref{thm_main2}. 
We thank I.~Domanov for pointing out that improved specific identifiability results
can be obtained by considering reshapings of the tensor.
The first and second author wish 
to thank the Simons Institute for the Theory of Computing in Berkeley, CA for generous support.

\section{The exceptional cases}\label{sec:newexample}
The following classical result shows that the values $n=2, 3, 5$, which appear
in Theorems \ref{thm_main} and \ref{thm_main2}, have a special role for elliptic
normal curves.

\begin{prop0}[Coble \cite{Cob}]\label{prop:n235}
Assume that there are only finitely many elliptic normal curves passing through
$k$ 
general points in $\Pj^n$.
Then, $n=2$, $3$, or $5$ and, correspondingly, $k=\frac{(n+1)^2}{n-1}$.
In these three cases, there is a unique elliptic normal curve in $\Pj^n$
passing through $\frac{(n+1)^2}{n-1}$ general points. 
\end{prop0}
\begin{proof}
Elliptic normal curves of degree $(n+1)$ in $\Pj^{n}$ depend on $(n+1)^2$ 
parameters, which is the dimension of the space of sections of the normal
bundle.
The passage of the curve through a point in $\Pj^{n}$ imposes $n-1$ 
conditions, which is the codimension of the curve.
Therefore, we may expect finitely many elliptic normal curves through $k$ 
\emph{general} points in $\Pj^{n}$ only if $k(n-1)=(n+1)^2$.
This implies that $(n-1)$ divides $(n+1)^2=(n-1)(n+3)+4$, hence $(n-1)$ divides 
$4$, which gives the values $n=2, 3, 5$. Moreover, $k=(n+1)^2/(n-1)$.

In case $n=2$ and $k=9$, the elliptic curve is a plane cubic, and it is unique. 

In case $n=3$ and $k=8$, an elliptic normal curve is a complete
intersection of two quadrics. Thus, if $\langle Q_1, Q_2\rangle$ is the pencil
of quadrics through $8$ general points $p_1,\ldots, p_8$, then $C=Q_1\cap Q_2$
is the unique elliptic normal curve through the $p_i$'s. 

In case $n=5$ and $k=9$, the existence and the uniqueness of the curve was 
found by Coble \cite[Theorem 19]{Cob} by applying a Gale transform---see
\cite{EP} for a nice review---and reducing to the case $n=2$ and $k=9$; a modern
treatment was given by Dolgachev \cite[Theorem 5.2]{Dolg}.
\end{proof}

Next, we analyze the case of cubic hypersurfaces in $\Pj^5$ that are singular 
at $9$ general points.

\begin{prop0}[Veneroni {\cite[Section 1]{Veneroni1905}}, Coble {\cite[p.
16]{Cob}}, Room {\cite[Sections 9--22]{Room1938}}, Fisher {\cite[Lemma
2.9]{Fisher}}]\label{prop:veneroni}
We have the following two results.
\begin{enumerate}[(i)]
\item The $2$-minors of a $3\times 3$ matrix with linear entries on $\Pj^5$ 
define a (sextic) elliptic normal curve in $\Pj^5$. 
\item If $\Var{C}$ is a (sextic) elliptic normal curve in $\Pj^5$, then 
 the variety of secant lines
$\Sec{2}{\Var{C}}$ is a complete intersection of two cubic 
hypersurfaces on $\Pj^5$, each 
one being the determinant of a $3\times 3$ matrix with linear entries on
$\Pj^5$.
\end{enumerate}
\end{prop0}
\begin{proof} Part (i) is well known: The curve is obtained by cutting the
Segre variety $\Pj^2\times\Pj^2\subset \Pj^9$, i.e., the variety of $3\times 3$ matrices of 
rank $1$, with a linear space $\Pj^5$. Claim (ii) follows by \cite[Lemma 2.9]{Fisher}.
\end{proof}

\begin{thm0}\label{thm:elliptic9} Let $p_1, \ldots, p_9$ be general points 
in $\Pj^5$. Let $\Var{C}$ be the elliptic normal sextic curve through these
points. A cubic that is singular at $p_1, \ldots, p_9$ contains  $\Sec{2}{\Var{C}}$ 
and is singular on $\Var{C}$.
\end{thm0}
\begin{proof} By \refprop{prop:veneroni}, in the pencil of cubics containing $\Sec{2}{\Var{C}}$, the 
general element is singular along $\Var{C}$.
This pencil fills the space of cubics that are singular at $p_1, \ldots, p_9$,
which is $2$-dimensional by the Alexander--Hirschowitz theorem \cite{AH1995}.
\end{proof}

These observations lead to a different proof of the third exceptional case in \refthm{thm_main}.

\begin{prop0}[Ranestad--Voisin \cite{Ranestad2013}]
The general tensor in $S^3\C^6$ of rank $9$ has exactly two Waring 
decompositions as sum of $9$ powers of linear forms.
\end{prop0}
\begin{proof}
In the language of \cite{ChCi2006}, we have to prove that the
secant order of $\Sec{9}{v_3(\Pj^5)}$ is $2$. By \cite[Theorem 2.4]
{ChCi2006}, this is equal to the secant order of the 
 $9$-contact locus $\Var{C}$, which corresponds to the third Veronese embedding 
 of an elliptic normal sextic curve in $\Pj^5$, by
Theorem \ref{thm:elliptic9}.
Thus, $\Var{C}$ is an elliptic curve of degree $18$ 
in $\Pj^{17}$, whose secant order is $2$ by \cite[Proposition 5.2]{ChCi2006}.
\end{proof}

\begin{coro0}\label{coro:veneroni}
Let $n=2$, $3$, or $5$.
The unique elliptic normal curve that is mentioned in \refprop{prop:n235}, 
which passes through $\tfrac{(n+1)^2}{n-1}$ general points in $\Pj^n$,
can be constructed as the singular locus of a general hypersurface of 
degree $\tfrac{2(n+1)}{n-1}$ that is singular in the $\tfrac{(n+1)^2}{n-1}$ points.
\end{coro0}
\begin{proof} The cases $n=2, 3$ have already been considered in the proof of
\refprop{prop:n235}. The case $n=5$ follows from \refthm{thm:elliptic9}.
\end{proof}

\begin{rem0}\label{rem:genesis} 
It was only at the completion of this manuscript that we became aware of
Ranestad and Voisin's proof of the third exceptional case that appears in 
Theorems \ref{thm_main} and \ref{thm_main2}. 
Our initial motivation for studying this problem arose because the third case
was unexpectedly---in our minds---suggested by a computational analysis 
performed by the third author, who ran the algorithm that we present in \refsec{sec:algorithm}, 
for hypersurfaces of degree $d$ in $\Pj^n$ singular at the maximal number of 
random points, i.e., $r = r_{d,n}-1$ with $r_{d,n}$ as in \refeqn{eqn_exp_rank}, 
for all reasonably small values of $d$, $n$.
It took a while to realize what happened, because this third case was missing
in \cite[Theorem 6.1.2]{LandsbergOttaviani2013}. Actually, Theorem 6.1.2 
of \cite{LandsbergOttaviani2013} only intended to collect previous results by 
Ballico \cite{BL}, Ciliberto and Chiantini \cite{ChCi2001}, and 
Mella \cite{Mella2006, Mella2009}, which are individually correct. 
The second author takes the responsibility to have first overlooked the assumption 
$d\ge 4$ in summarizing and reporting the results of \cite{Mella2006, BL}. 
For $d = 3$, Theorem \ref{thm_main} was known with the additional assumption
$r < r_{3,n} -\frac{n+2}{3}+1$; see \cite[Theorem 4.1]{Mella2006}. 
From the theoretical proof that we present in \refsec{sec:cubics}, 
we can conclude that the third case was the last exception.
Therefore, Theorem 6.1.2 in \cite{LandsbergOttaviani2013} remains true if the 
third case $(k,d,n)=(9,5,3)$ is added to the list of exceptions. 
Exactly the same remark applies to the formulation of Theorem 2.3
in \cite{OeOtt} and Theorem 12.3.4.3 in \cite{Land_book}. 
We informed the coauthors of \cite{LandsbergOttaviani2013, OeOtt} of the 
problem, and they accepted the above conclusion.
\end{rem0}

\section{Cubics singular at the maximum number of points.}\label{sec:cubics}
We turn our attention to the proof of \refthm{thm_main} in the case of cubics,
i.e., $d=3$. Given $n$, we define 
\begin{align*}
k_n = \Bigl\lceil\frac{\binom{n+3}{3}}{n+1} \Bigr\rceil = \Bigl\lceil
\frac{(n+3)(n+2)}{6}\Bigr\rceil;
\end{align*}
it is the generic rank for cubic polynomials for $n\neq 4$. In other words,
a cubic polynomial on $\Pj^n$ singular at $k_n$ general points vanishes 
identically for $n\neq 4$ \cite{AH1995}. Some elementary algebra shows that 
$k_n=\frac{(n+3)(n+2)}{6}$ if $n \not\equiv 2 \mod 3$, while 
$k_n=\frac{(n+3)(n+2)}{6}+\frac{2}{3}=\frac{(n+4)(n+1)}{6}+1$ if $n\equiv 2 \mod 3$. 

For the sake of future reference, let us state explicitly the following consequence 
of the Alexander--Hirschowitz theorem \cite{AH1995}.

\begin{thm0}[Alexander--Hirschowitz \cite{AH1995}] \label{AHcubics}
The space of cubic hypersurfaces on $\Pj^n$ that are singular at $k_n-1$ 
general points has dimension
\begin{enumerate}[(i)]
\item $n+1$ if $n \not\equiv 2 \mod 3$, or
\item $\tfrac{n+1}{3}$ if $n \equiv 2 \mod 3$.
\end{enumerate}
In addition, the space of cubic hypersurfaces on $\Pj^n$, $n\ne4$, that are singular at 
$k_n$ general points is empty.
\end{thm0}

To complete the proof of the Theorems \ref{thm_main} and \ref{thm_main2}, it remains to show the 
following result, which refines Theorem \ref{AHcubics}.

\begin{thm0}\label{cubiche}
The space of cubic hypersurfaces in $\Pj^n$ that are singular at $k_n-1$ 
general points has dimension
\begin{enumerate}[(i)]
\item $n+1$ if $n \not\equiv 2 \mod 3$, or
\item $\tfrac{n+1}{3}$ if $n\equiv 2 \mod 3$,
\end{enumerate}
and, in addition, its common singular locus consists only of these $k_n-1$ 
points, provided that $n\neq 5$.
\end{thm0}

In order to prove \refthm{cubiche}, we may assume $n\ge 6$, since
the cases with $n\le 4$ (as well as the case of cubics in $\Pj^5$ singular at $8$ points)
can be checked separately using the approach described in \refsec{sec:algorithm}. 
The outline of our proof of \refthm{cubiche} is as follows. In 
\refsec{sec:codim3}, we will prove case (i) by induction on subspaces
of codimension $3$, adopting an approach that is mainly inspired 
by \cite[Section 5]{Brambilla2008}, where an alternative proof of Theorem \ref{AHcubics} was given. 
To prove case (ii), the aforementioned 
technique needs a modification. We will construct an inductive proof 
on subspaces of codimension $3$ and $4$; in the inductive step, we will rely,
additionally, on the argument of case (i). This strategy will be presented in 
\refsec{sec:codim4}.

In the rest of this section, if $S$ is a set of simple points in $\Pj^n$ 
and $P \subset \Pj^n$ is a linear subspace, we denote with $I_{S,P}(d)$ 
the space of degree $d$ polynomials in $P$ vanishing at all of the points in $S$. 
Moreover, if $\mathbf{X}$ is a a set of double (singular) points, we 
denote by $I_{\mathbf{X} \cup S,P}(d)$  the space of degree $d$ polynomials 
in $P$ vanishing on all of the points in $S \cup \mathbf{X}$ and whose 
derivatives vanish on all of the points in $\mathbf{X}$. The notation 
$L = (x_{i} \ldots x_{i+d})$ denotes the subspace of codimension $d+1$ whose ideal 
is $\langle x_i, x_{i+1}, \ldots, x_{i+d} \rangle$.

\subsection{Proof of Theorem \ref{cubiche} (i) by induction on codimension $3$.}\label{sec:codim3}
We start by proving three auxiliary results.

\begin{prop0}\label{proprep}
Let $n\ge 6$, and let $L, M, N \subset\Pj^ n$ be general subspaces of 
codimension $3$. 
Let $l_i$, respectively $m_i$, with $i=1, 2, 3$ be three general 
points on $L$, respectively $M$.
Let $n_i$ with $i=1, 2$ be two general points on $N$. 
Then, the space of cubic hypersurfaces in $\Pj^n$ that contain 
$L \cup M \cup N$ and that are singular at the eight points 
$\mathbf{X} = \{l_1, l_2, l_3, m_1, m_2, m_3, n_1, n_2\}$ has dimension $3$.
Furthermore, the common singular locus is contained in $L\cup M\cup N$.
\end{prop0}
\begin{proof}
The base cases $n=6$, $7$ and $8$ can be proved with the Macaulay2 script
\texttt{generic-identifiability.m2} that is provided as an ancillary file 
to the arXiv version of this paper. Using this software, we may compute 
the following  dimensions:
\begin{align*}
\dim I_{\mathbf{X}\cup L\cup N\cup M, \Pj^6}(2) &= 0, &\dim I_{\mathbf{X}
\cup L\cup N\cup M, \Pj^6}(3)&= 3, \;\text{and} \\
\dim I_{\mathbf{X}\cup L\cup N\cup M, \Pj^7}(2) &= 0, &\dim I_{\mathbf{X}
\cup L\cup N\cup M, \Pj^7}(3)&= 3,
\end{align*}
so that the claim on the codimension follows.
The code also proves the statement about the singular locus.

For $n\ge 9$, the statement follows by induction on $n$.
Indeed, we may choose coordinates such that $L=(x_0\ldots x_2)$, 
$M=(x_3\ldots x_5)$, $N=(x_6\ldots x_8)$. In this setting it is clear that there are no quadrics 
that contain $L\cup M\cup N$, and moreover every cubic containing $L\cup M\cup N$
is a cone with vertex in $L\cap M\cap N$.
Thus, for a general hyperplane $H\subset \Pj^n$, 
the Castelnuovo sequence (see \cite[Equation (1)]{Brambilla2008}) induces 
an inclusion
\begin{align*}
0 
\rig{} I_{L\cup M\cup N,\Pj^n}(3)
\rig{} I_{(L\cup M\cup N)\cap H,H}(3).
\end{align*}
Hence, if we specialize the eight points to the hyperplane $H$, we get
an inclusion
\begin{align*}
0 
\rig{} I_{\mathbf{X} \cup L \cup M \cup N, \Pj^n}(3) 
\rig{} I_{(\mathbf{X}\cup L\cup M\cup N)\cap H,H}(3).
\end{align*} 
Then, our statement follows by induction. The singular locus is a cone with vertex $L\cap
M\cap N$
over the singular locus of the base case $n=8$.
\end{proof}

\begin{rem0}
Following the output of the software for the case $n=8$, we can guess
the common singular locus of the cubic hypersurfaces in 
$I_{\mathbf{X} \cup L \cup M \cup N, \Pj^n}(3) $, for $n\geq 8$. 
It turns out that in some examples---but we believe in general---it is 
given by the union of the three linear subspaces
$L\cap M$, $L\cap N$, $M\cap N$ and by $8$ linear subspaces of codimension 
$7$, each containing one of the $8$ points, and
three of them contained in $L$, three of them contained in $M$, and two 
of them contained in $N$.
\end{rem0}

\begin{prop0}\label{proprep2}
Let $n\ge 5$, and let ${L, M}\subset\Pj^ n$ be subspaces of codimension three.
Let $l_i$, respectively $m_i$, with $i=1, \ldots, n-2$ be general points on
${L}$,
respectively ${M}$. Let $p_1, p_2 \in \Pj^n$ be general points. Then, 
the space of cubic hypersurfaces in $\Pj^n$ 
containing $L\cup M$ and singular along the set of $2n-2$ points 
$\mathbf{X} = \{l_1, l_2, \ldots, l_{n-2}, m_1, m_2, \ldots, m_{n-2}, 
p_1, p_2 \}$ has dimension $n+1$. 
Its common singular locus contains the linear space $L\cap M$ 
and is $0$-dimensional at the points $p_1$ and $p_2$.
\end{prop0}
\begin{proof}
The base cases $n=5$, $6$, and $7$ can be proved with the Macaulay2 
script \texttt{generic-identifiability.m2}.
Running the software, we find the following dimensions
\begin{align*}
 \dim I_{\mathbf{X} \cup L \cup M,\Pj^5} (2) &= 0, &\dim I_{\mathbf{X} \cup L 
 \cup M,\Pj^5} (3) &= 6, \\
 \dim I_{\mathbf{X} \cup L \cup M,\Pj^6} (2) &= 0, &\dim I_{\mathbf{X} \cup L 
 \cup M,\Pj^6} (3) &= 7, \text{ and }\\
 \dim I_{\mathbf{X} \cup L \cup M,\Pj^7} (2) &= 0, &\dim I_{\mathbf{X} \cup L 
 \cup M,\Pj^7} (3) &= 8.
\end{align*}
These values indeed correspond to the claimed dimensions.

For $n\ge 8$, the statement follows by induction from $n-3$ to $n$.
Indeed, given a third general subspace ${N}$ of codimension $3$, we get 
the exact sequence
\begin{align*}
0
\rig{} I_{L\cup M\cup N,\Pj^n}(3)
\rig{} I_{L\cup M,\Pj^n}(3) 
\rig{} I_{(L\cup M)\cap N,N}(3),
\end{align*}
where the dimensions of the three spaces in the sequence are respectively $27$,
$9(n-1)$, and $9(n-4)$. 
Let us specialize $n-5$ of the points $l_i \in L$ to $L\cap N$, $n-5$ 
of the points $m_i \in M$ to $M\cap N$, 
and the two points $p_1,p_2$ to $N$. 
Then, we obtain a sequence
\begin{align*}
0 
\rig{} I_{\mathbf{X}\cup L\cup M\cup N,\Pj^n}(3) 
\rig{} I_{\mathbf{X}\cup L\cup M,\Pj^n}(3) 
\rig{} I_{(\mathbf{X}\cup L\cup M)\cap N,N} (3),
\end{align*}
where the trace $\left(\mathbf{X}\cup{ L}\cup{ M}\right)\cap { N}$
 satisfies the assumptions on $N=\Pj^{{n-3}}$, so that we can apply 
 the induction.
Notice that the residual (left) space satisfies the hypotheses 
of \propref{proprep} and has dimension $3$.
Since the common singular locus of the cubics containing $L\cup M$ and 
singular at $\mathbf X$ must be contained in the common singular locus 
of the leftmost $3$-dimensional space, 
it follows by \propref{proprep} that its components through $p_1$ and $p_2$ 
must be contained in $N$.
After the degeneration, the space of cubics $I_{\mathbf{X}\cup L\cup 
M,\Pj^n}(3)$ still has dimension at most
$3+(n-2)=n+1$ by induction. Hence, by semicontinuity it follows that 
its dimension is indeed equal to $n+1$. 
The common singular locus cannot be positive dimensional at points $p_1$ 
and $p_2$, because otherwise it should be of positive dimension in the trace 
(right space), where by induction we know that it is $0$-dimensional.
\end{proof}

\begin{prop0}\label{proprep3}
Let $n\ge 6$, and let ${L}\subset\Pj^n$ be a subspace of codimension $3$.
If $n \not\equiv 2 \mod 3$, then the space of cubic hypersurfaces in 
$\Pj^n$ that contain $L$ and 
that are singular at $\frac{n(n-1)}{6}$ general points $l_i \in L$ and at 
$n$ general points
$p_i \in \Pj^n$ has dimension $n+1$. Moreover, its common singular 
locus is $0$-dimensional at 
the $n$ points $p_i$.
\end{prop0}
\begin{proof}
The statement can be checked for $n=6$, $7$ with the Macaulay2 script 
\texttt{generic-identifiability.m2.}

Let $n\ge 9$ and $n \not\equiv 2 \mod 3$. Consider the sequence
\begin{align*}
0
\rig{} I_{L\cup M,\Pj^n}(3) 
\rig{} I_{L,\Pj^n}(3)
\rig{} I_{L \cap M,M}(3)
,
\end{align*}
where $M$ is a general subspace of codimension $3$. Denoting by $\mathbf{X}$
the union of the double points supported at the points $l_i$ and $p_i$, we get
\begin{align*}
0
\rig{} I_{\mathbf{X}\cup L\cup M,\Pj^n}(3)
\rig{} I_{\mathbf{X}\cup L,\Pj^n}(3)
\rig{} I_{(\mathbf{X}\cup L)\cap M,M}(3).
\end{align*}
We specialize $\frac{(n-3)(n-4)}{6}$ of the points $l_i$ to $L\cap M$ and 
$n-2$ of the points $p_i$ to $M$. We can assume that at least one point 
$p_i$ that is a contained in a positive dimensional component of the singular 
locus is not specialized. Thus, we left $n-2$ general points on $L$ and $2$ 
general points in $\Pj^n$. Let us note that we cannot apply 
induction from $n-3$ to $n$ 
to determine the dimension of the right space---contrary to the strategy 
that was employed in the proof of the foregoing propositions in this 
section---because then we would have 
to specialize $n-3$ of the points $p_i$ to $M$, hereby losing control over the 
singular locus. Instead, we note that we can immediately use Proposition 
5.4 of \cite{Brambilla2008} (in $\Pj^{n-3}$) on the trace (right space); 
it turns out to be empty. 
On the residual (left space), \refprop{proprep2} can be invoked, proving 
that it has dimension $n+1$. 
Now, if the singular locus would have a positive dimensional component, then, since 
the dimension of the space of cubics is constant 
along the specialization (it equals $n+1$), we would get 
a deformation of the singular locus, which should 
 be of positive dimension at every point. This, however,
contradicts \refprop{proprep2}, hereby concluding the proof.
\end{proof}

We are now ready to prove the first part of \refthm{cubiche}.
\begin{proof}[Proof of \refthm{cubiche}, case (i).]
We fix a linear subspace $L\subset{\Pj}^n$ of codimension $3$ and
consider the exact sequence
$$
0
\rig{} I_{L,\Pj^n}(3)
\rig{} S_{\Pj^ n}(3)
\rig{} S_L(3),
$$
where $S_{\Pj^ n}(3)$ is the space of cubic polynomials on $\Pj^n$
and the quotient space $S_L(3)$ is isomorphic to
the space of cubic polynomials on $L$.
Then, we specialize to $L$ as many points as possible in such a 
way that the trace with respect to $L$ imposes independent conditions on the 
cubics of $L$. To be precise, we have $k_n-1=\frac{(n+3)(n+2)}{6}-1$ double 
points and we specialize $k_{n-3}=\frac{n(n-1)}{6}$ of them to $L$, 
leaving $n$ points outside. 
Then, the result follows from Theorem 5.1 of \cite{Brambilla2008} on 
the trace (right space), which turns out to be empty, and by 
\refprop{proprep3} on the residual (left space), which has dimension $n+1$. 
If the contact locus has positive dimension, then, since
the dimension of the space of cubics is constant and equal to $(n+1)$ in the degeneration, 
we would get a deformation of the singular locus with a positive dimension at every point,
contradicting \refprop{proprep3} and concluding the proof.
\end{proof}

\subsection{Proof of Theorem \ref{cubiche} (ii) by induction on codimension $3$ and $4$.}
\label{sec:codim4}
For proving the second case in \refthm{cubiche}, we need to introduce 
several other auxiliary results on configurations that involve subspaces 
of codimension three and four. These configurations are covered in 
Propositions \ref{proprep:codim433} through \ref{proprep3:codim4}.

\subsubsection{Codimension $4$, $3$, $3$}
\begin{prop0} \label{proprep:codim433}
Let $n\ge 6$, and let $L, M, N \subset\Pj^ n$ be general subspaces of
codimension $4$, $3$, and $3$, respectively. Let $l_1, l_2, l_3$ be 
general points on $L$.
Let $m_i$, respectively, $n_i$ with $i=1, \ldots, 4$ be four general 
points on $M$, respectively $N$. 
Then, the space of cubic hypersurfaces in $\Pj^n$ that contain $L\cup M\cup N$ 
and are singular at the $11$ points $\mathbf{X} = \{l_1,l_2,l_3, 
m_1,\ldots,m_4, n_1,\ldots,n_4\}$ is empty.
\end{prop0}
\begin{proof} The proof is similar to the proof of \refprop{proprep}.
The Macaulay2 code  proves the base cases $n=6$, $7$, $8$ and $9$.
For $n\ge 9$, we may choose coordinates such that $L=(x_0\ldots x_3)$, $M=(x_4\ldots x_6)$, and $N=(x_7\ldots x_9)$. 
Then, the statement follows by induction on $n$.
Indeed, as in the proof of \refprop{proprep}, the space $I_{{ L}\cup{ M}\cup{ N},{\Pj^ n}}(2)$ is empty,
thus for a general hyperplane $H\subset \Pj^n$, the Castelnuovo 
sequence induces an embedding
\begin{align*}
0
\rig{} I_{L\cup M\cup N,\Pj^n}(3)
\rig{} I_{(L \cup M\cup N)\cap H,H}(3),
\end{align*}
and, moreover, every cubic in the left space is a cone with vertex at $L\cap M\cap N$.
Hence, by specializing the $11$ points to the hyperplane $H$, we get:
\begin{align*}
0 
\rig{} I_{\mathbf{X}\cup L\cup M\cup N,\Pj^n}(3)
\rig{}I_{(\mathbf{X}\cup L\cup M\cup N)\cap H,H}(3).
\end{align*}
Then, the statement follows by induction.
\end{proof}

\begin{prop0}\label{proprep2:codim433}
Let $n\ge 7$, let $n \equiv 1 \mod 3$, and let $L, M \subset \Pj^n$ be 
subspaces of codimension $4$ and $3$, respectively.
Let $l_i$ with $i=1, \ldots, n-3$ be general points on ${L}$.
Let $m_i$ with $i=1, \ldots, \frac{4n-10}{3}$ be general points on ${M}$.
Then, the space of cubic hypersurfaces in $\Pj^n$ that contain $L \cup M$ 
and are singular at all the $l_i$'s and $m_i$'s and at four general 
points $p_1,p_2,p_3,p_4 \in \Pj^n$, is empty.
\end{prop0}
\begin{proof}
The Macaulay2 script proves the base case $n=7$.

For $n = 3 k + 1$ with $k \ge 3$, the statement follows by induction 
from $n-3$ to $n$. Indeed, given a third general subspace ${N}$ of 
codimension $3$, we get the exact sequence
\begin{align*}
0
\rig{} I_{L\cup M\cup N,\Pj^n}(3)
\rig{} I_{L\cup M,\Pj^n}(3)
\rig{} I_{(L \cup M)\cap N,N}(3),
\end{align*}
where the dimensions of the three spaces in the sequence are respectively 
$36$,  $12n-18$ and $12n-54$.
Let $\mathbf{X}$ denote the union of the double points supported at the 
$p_i$'s, $l_i$'s and $m_i$'s.
Assume that we specialize $n-6$ of the points $l_i \in L$ to $L\cap N$,
$\frac{4n-22}{3}$ of the points $m_i \in M$ to $M\cap N$, and the four 
points $p_1,\ldots, p_4$ to $N$. Then, we obtain a sequence
\begin{align*}
0
\rig{} I_{\mathbf{X}\cup L\cup M\cup N,\Pj^n}(3)
\rig{} I_{\mathbf{X}\cup L\cup M,\Pj^n}(3)
\rig{} I_{(\mathbf{X}\cup L\cup M)\cap N,N}(3),
\end{align*}
where the trace $(\mathbf{X}\cup L\cup M)\cap N$ satisfies 
the assumptions on $N=\Pj^{{n-3}}$, so that we can apply induction.
Then, we may conclude, as the residual (left space) satisfies the 
hypotheses of \propref{proprep:codim433}, and, consequently, it is empty.
\end{proof}

\begin{prop0}\label{proprep3:codim433}
Let $n\ge 7$, $n \equiv 1 \mod 3$, and ${L}\subset\Pj^n$ be a subspace 
of codimension four. 
Then, the space of cubic hypersurfaces in $\Pj^n$ that are singular at
$k_{n-4}=\frac{(n-1)(n-2)}{6}$ general points $l_i$ on ${L}$ 
(and, thus, contain ${L}$, by \refthm{AHcubics}) 
and at $\frac{4n+2}{3}$ general points $p_i\in\Pj^n$
is empty.
\end{prop0}
\begin{proof}
The Macaulay2 script \texttt{generic-identifiability.m2} proves the case $n=7$.

For $n = 3k + 1$ with $k \ge 3$, the statement follows by the sequence
\begin{align*}
0
\rig{} I_{L \cup M,\Pj^n}(3)
\rig{} I_{L,\Pj^n}(3)
\rig{} I_{L \cap M,M}(3),
\end{align*}
where $M$ is a general subspace of codimension $3$. 
If we denote by $\mathbf{X}$ the union of the double points supported 
at the points $l_i$ and $p_i$, then we get the sequence
\begin{align*}
0
\rig{} I_{\mathbf{X} \cup L \cup M,\Pj^n}(3)
\rig{} I_{\mathbf{X} \cup L,\Pj^n}(3)
\rig{} I_{(\mathbf{X} \cup L) \cap M,M}(3).
\end{align*}
Then, we specialize $k_{n-7}$ of the points $l_i$ to $L\cap M$ and 
$\frac{4n-10}{3}$ of the points $p_i$ to $M$.
The trace (right space) contains exactly $k_{n-3}$ double points and 
turns out to be empty by induction.
Thus, there remain $n-3$ general points on $L$ and $4$ general points 
on $\Pj^n$; we can then use \refprop{proprep2:codim433} on the residual 
(left space) to conclude.
\end{proof}

\subsubsection{Codimension 4, 4, 3}
\label{subsection:443}

\begin{prop0}\label{proprep:codim4}
Let $n\ge 8$, and let $L, M, N \subset\Pj^ n$ be general subspaces of
codimension respectively $4$, $4$, and $3$. 
Let $l_i$, respectively $m_i$, with $i=1,\ldots, 4$ be general points on $L$,
respectively $M$.
Finally, let $n_i$ with $i=1, \ldots, 5$ be general points on $N$. 
Then, the space of cubic hypersurfaces in $\Pj^n$ that contain $L \cup M 
\cup N$ and that are singular at the $13$ points $\mathbf{X} = \{ l_1, 
\ldots, l_4, m_1, \ldots, m_4, n_1, \ldots, n_5\}$ 
has dimension $1$. In other words, there is a unique cubic hypersurface $W$
through $L \cup M \cup N$ and singular at $\mathbf{X}$. 
Furthermore, the singular locus of $W$ is contained in $L\cup M\cup N$.
\end{prop0}
\begin{proof}
The proof is similar to the proof of \refprop{proprep}.
The Macaulay2 code \texttt{generic-identifiability.m2} proves the base cases $n=8$, $9$, and $10$.

For $n\ge 11$, we may choose coordinates such that $L=(x_0\ldots x_3)$, $M=(x_4\ldots x_7)$, 
and $N=(x_8\ldots x_{10})$; then, the statement follows by induction on $n$. Indeed, as in the proofs of 
\refprop{proprep} and \refprop{proprep:codim433}, we let $H\subset \Pj^n$ be a general hyperplane, 
so that the Castelnuovo sequence induces the inclusion
$$
0
\rig{} I_{L \cup M \cup N,\Pj^n}(3)
\rig{} I_{(L \cup M \cup N) \cap H,H}(3),
$$
because the space $I_{L \cup M \cup N,\Pj^n}(2)$ is empty. 
Hence, by specializing the $13$ points on the hyperplane $H$, we get an
exact sequence:
$$
0 
\rig{} I_{\mathbf{X} \cup L\cup M \cup N,\Pj^n}(3)
\rig{} I_{(\mathbf{X} \cup L \cup M \cup N)\cap H,H}(3).
$$
Now the statement follows by induction.
\end{proof}

\begin{prop0}\label{proprep2:codim4}
Let $n\ge 8$,  $n \equiv 2 \mod 3$, and ${L, M}\subset\Pj^n$ be subspaces 
of codimension $4$.
Let $l_i$ and $m_i$, where $i=1, \ldots, \frac{4n-14}{3}$, be general 
points on ${L}$ and ${M}$, respectively. 
Then, the space of cubic hypersurfaces in $\Pj^n$ that contain $L \cup M$ 
and are singular
at the $\frac{8n-28}{3}$ points $l_i, m_i$, $i=1,  \ldots, \frac{4n-14}{3}$, 
and at an additional set of five 
general points $p_i \in \Pj^n$, $i=1,\ldots, 5$, has dimension $\frac{n+1}{3}$. 
Furthermore, its common singular locus, which contains the linear space 
$L\cap M$,  is $0$-dimensional at each of the points $p_1, \ldots, p_5$.
\end{prop0}
\begin{proof}
The case $n=8$ is handled in the \texttt{generic-identifiability.m2} Macaulay2 script.

For $n = 3k + 2$ with $k\ge3$, the statement follows by induction on $k$. 
Given a third general subspace $N$ of codimension $3$, we get the exact sequence
$$
0
\rig{} I_{L \cup M \cup N,\Pj^n}(3)
\rig{} I_{L \cup M,\Pj^n}(3)
\rig{} I_{(L \cup M) \cap N,N}(3),
$$
where the dimensions of the three spaces in the sequence are 
respectively $48$,  $16(n-2)$ and $16(n-5)$.

Let $\mathbf{X}$ denote the union of the double points supported at 
$p_1,\ldots, p_5$, $l_i$ and $m_i$ with $i=1,\ldots,\frac{4n-14}{3}$. 
Then, we specialize $\frac{4n-26}{3}$ of the points $l_i \in L$ to $L\cap N$, 
$\frac{4n-26}{3}$ of the points $m_i \in {M}$ to $M\cap N$, 
and the points $p_1,\ldots, p_5$ to $N$.
We thus obtain a sequence
$$
0 
\rig{} I_{\mathbf{X} \cup L \cup M \cup N,\Pj^n}(3)
\rig{} I_{\mathbf{X} \cup L \cup M,\Pj^n}(3)
\rig{} I_{(\mathbf{X} \cup L \cup M) \cap N,N}(3),
$$
where the trace $\left(\mathbf{X}\cup{ L}\cup{ M}\right)\cap { N}$
satisfies the assumptions on $N=\Pj^{{n-3}}$, so that we can apply induction.
Then, the residual (left space) satisfies the hypotheses 
of \propref{proprep:codim4} and has dimension one. Moreover,
the common singular locus has to be contained in the common singular 
locus of the left $1$-dimensional space.
After the degeneration, the space $I_{\mathbf{X} \cup L \cup M,\Pj^n}(3)$ 
still has dimension less than 
or equal to $1+\frac{n-2}{3}=\frac{n+1}{3}$, by induction, and, therefore, 
its dimension equals $\frac{n+1}{3}$, by semicontinuity. 
The common singular locus cannot be positive dimensional at the points 
$p_1,\ldots, p_5$, 
because otherwise it should be positive dimensional in the trace (right space),
while we know that it is $0$-dimensional there by induction.
\end{proof}

\begin{prop0}\label{proprep3:codim4}
Let $n\ge 8$, $n \equiv 2 \mod 3$, and ${L}\subset\Pj^n$ be a subspace 
of codimension $4$. 
Then, the space of cubic hypersurfaces in $\Pj^n$ that contain $L$ 
and are singular at 
$k_{n-4}=\frac{(n-1)(n-2)}{6}$ general points $l_i \in L$ and at 
$\frac{4n+1}{3}$ general points $p_i \in \Pj^n$
has dimension $\frac{n+1}{3}$. 
Furthermore, its singular locus is of dimension $0$ at all of the points $p_i$.
\end{prop0}
\begin{proof}
The statement follows by the sequence
$$
0
\rig{} I_{L \cup M,\Pj^n}(3)
\rig{} I_{L,\Pj^n}(3)
\rig{} I_{L \cap M,M}(3),
$$
where $M$ is a general subspace of codimension $4$. Denoting by $\mathbf{X}$
the union of the double points supported at the points $l_i$'s and $p_i$'s, 
we get
$$
0
\rig{} I_{\mathbf{X} \cup L \cup M,\Pj^n}(3)
\rig{} I_{\mathbf{X} \cup L,\Pj^n}(3)
\rig{} I_{(\mathbf{X} \cup L) \cap M,M}(3).
$$

Suppose that the singular locus would be of positive dimension at one of the $p_i$'s, 
say, at $q = p_j$.
Then, we specialize $k_{n-8}$ of the points $l_i$ to $L\cap M$ and 
$\frac{4n-14}{3}$ of the points $\{p_i\}_{i\ne j}$ to $M$.
The trace (right space) contains exactly $k_{n-4}$ double points and 
is empty because of \refprop{proprep3:codim433}.
There remain $\frac{4n-14}{3}$ general points on $L$ and $5$ general points---one 
of which is $q$---on $\Pj^n$. We can use \refprop{proprep2:codim4} on the 
residual (left space), which has dimension $\frac{n+1}{3}$. 
By assumption, the singular locus has positive dimension 
at $q$; however, the dimension of the space of cubics is constant and equal to 
$\frac{n+1}{3}$ through the degeneration, so that we have 
a deformation of the singular locus, which is 
of positive dimension at every point. In particular, the singular locus will
be of positive dimension at all the points that were not specialized to $M$, hereby 
contradicting \refprop{proprep2:codim4}. We conclude that our initial assumption 
must have been false, so that no general points $p_j$ can exist where the 
singular locus is of positive dimension.
\end{proof}

\begin{proof}[Proof of Theorem \ref{cubiche}, part (ii).]
We fix a codimension four linear subspace $L\subset \Pj^n$ and we
use the exact sequence
$$
0
\rig{} I_{{L},{\Pj^ n}}(3)
\rig{} S_{\Pj^ n}(3)
\rig{} S_L(3),
$$
where, as above, $S_{\Pj^ n}(3)$ is the space of cubic 
polynomials on $\Pj^n$ and the quotient space $S_L(3)$ is isomorphic to
the space of cubic polynomials on $L$.
We specialize $k_{n-4}$ points on $L$, leaving
$\frac{4n+1}{3}$ points outside. Then, the result follows 
from Theorem 5.1 of \cite{Brambilla2008} on the trace
(right space),
which turns out to be empty and
by \refprop{proprep3:codim4} on the residual (left space).
If the contact locus would have a positive dimension, then, since the 
dimension of the space of cubics is constant and equal to $\frac{n+1}{3}$ in the degeneration, 
we would get a deformation of the singular locus, which should be of 
positive dimension at every point; however, this contradicts \refprop{proprep3:codim4}, 
hereby concluding the proof.
\end{proof}

\section{Dual varieties to the relevant secant varieties}\label{sec:dual}
Denote by $\Tang{x}{\Var{X}}$ the tangent space to the projective 
variety $\Var{X} \subset \Pj^n$ at the point $x \in \Var{X}$.
Following the notation of \cite{ChCi2006} 
we say that $\Var{X}$ is \emph{not $k$-weakly defective} if the general 
hyperplane $H$ containing the linear span of the tangent spaces at $k$ general 
points $x_1,\ldots, x_k \in \Var{X}$, i.e., 
$\langle \Tang{x_1}{\Var{X}}, \ldots, \Tang{x_k}{\Var{X}} \rangle \subset H$, is tangent to $\Var{X}$ 
only at finitely many points. This is equivalent with saying that 
the $k$-contact locus with respect to $x_1,\ldots, x_k$
and $H$ is zero-dimensional.

For any projective variety $\Var{X}$, we will denote by $\Var{X}^\vee$ 
the dual variety to $\Var{X}$. Note that the dual of the secant variety 
$\sigma_k(v_d(\Pj^n))^{\vee}$ contains the points corresponding to 
hypersurfaces of degree $d$ in $\Pj^n$ with $k$ general singular points,
and it has codimension $\ge k$, where $k$ is the expected value 
for the codimension.

\begin{prop0} Let $\Var{X}\subset\Pj^N$ and let $\sigma_k(\Var{X})$ be the $k$-secant variety 
of $\Var{X}$. Then, the following are equivalent:
\begin{enumerate}[(i)]
\item the general hyperplane $H$ containing $\langle \Tang{x_1}{\Var{X}}, \ldots, 
\Tang{x_k}{\Var{X}}\rangle$ for general $x_1,\ldots, x_k$ is tangent 
to $\Var{X}$ only at $x_1,\ldots, x_k$, i.e., the $k$-contact locus with respect 
to $x_1,\ldots, x_k$ and $H$ consists exactly of the points $x_1,\ldots, x_k$,
\item $\Var{X}$ is not $k$-weakly defective, and
\item $\dim \left[\sigma_k(\Var{X})\right]^{\vee}=N-k$, that is a general 
hyperplane tangent to $\sigma_k(\Var{X})$ is tangent along a linear space 
of projective dimension $k-1$.
\end{enumerate}
\end{prop0}

\begin{proof} (i) $\Longleftrightarrow$ (ii) follows from \cite[Theorem 1.4]
{ChCi2001}. (ii) $\Longleftrightarrow$ (iii) follows from Terracini's Lemma.
\end{proof}

It is interesting to describe the dual varieties of $\sigma_k(v_d(\Pj^n))$ 
in the exceptional cases of Theorems \ref{thm_main} and \ref{thm_main2}.
They have dimension smaller than expected.

\begin{thm0}\label{thm:dualveronese}
The following dual varieties correspond to the exceptional cases appearing in Theorems \ref{thm_main} and \ref{thm_main2}.
\begin{enumerate}[(i)]
\item $\sigma_9(v_6(\Pj^2))^{\vee}$ contains the plane sextics which are 
double cubics. It has codimension $18$.

\item $\sigma_8(v_4(\Pj^3))^{\vee}$ contains the quartic surfaces which 
are reducible in a pair of quadrics. It has codimension $16$.

\item $\sigma_9(v_3(\Pj^5))^{\vee}$ contains the cubic $4$-folds which 
can be written as the determinant of a $3\times 3$
matrix with linear entries. It has codimension $18$.  
\end{enumerate}
\end{thm0}

To compute the dimension in third case, note that the Hilbert scheme 
of elliptic normal sextic curves in $\Pj ^5$
has dimension $36$. So the cubic hypersurfaces coming from this 
construction have dimension $37$, and $37+18=55={{8}\choose 3}-1$.

We remark that the defective Veronese varieties according to 
the classification of Alexander and Hirschowitz \cite{AH1995}  
(see \cite{Ot2009} for the equations of the defective secant varieties)
yield the following dual varieties
\begin{enumerate}[(i)]
\item $\sigma_{n(n+3)/2}(v_4(\Pj^n))^{\vee}$, for $n=2, 3, 4$, 
contains quartic hypersurfaces which are double quadrics.
It has codimension ${{n+2}\choose 3}\frac{n+7}{4}$.

\item $\sigma_7(v_3(\Pj^4))^{\vee}$ contains cubic $3$-folds which can 
be written as the determinant of a $3\times 3$
symmetric matrix with linear entries. It has codimension $13$.  
Indeed, it is birational to the Hilbert scheme  of quartic rational normal curves which has dimension
$21$.
\end{enumerate}

\section{Specific identifiability of symmetric tensors}
\label{sec:algorithm}

While the generic symmetric tensor of subgeneric rank is expected to admit a unique 
Waring decomposition, \emph{specific} tensors, whose Waring decomposition is assumed 
to be known, may admit
multiple decompositions. We proceed by presenting an approach for certifying 
specific identifiability of symmetric tensors of small rank by checking not tangential
weak defectivity of the $r$-secant variety of a Veronese variety in the given point. 
The strategy is an adaption of the algorithm from \cite{COV2014} to the setting 
of identifiability with respect to the Veronese variety $\Var{V} = v_d(\Pj^n)$. As such,
the presented condition will only be a sufficient condition; that is, if the criterion 
does not apply, then the outcome of the test is inconclusive. On the other hand, if the 
criterion applies, then the given input tensor is $r$-identifiable and of symmetric 
rank $r$. Throughout this section, it is assumed that we are handed a Waring decomposition 
\[
 p = p_1 + \cdots + p_r \in \Sec{r}{\Var{V}} \subset S^d \C^{n+1},
\]
wherein the point $p_i = \vect{a}_i^{\otimes d} \in \Var{V}$ is the degree $d$ Veronese embedding of  
the vector $\vect{a}_i \in \C^{n+1}$. In other words, we know the points $p_i$ appearing in the 
decomposition. The goal only consists of certifying that $p$ is $r$-identifiable---the 
decomposition(s) are not sought. To this end, the strategy in \cite{COV2014} suggests a 
two-step procedure: Prove that $p$ is a smooth point, and verify the Hessian criterion.
In principle, the method can be applied for all tensors of subgeneric rank, however, 
in practice the range of applicability of the algorithm in \refsec{sec_algorithm} is 
restricted by the lack of good techniques for certifying smoothness. For this reason, 
the results of \cite{DDL2013} may apply in a wider range than the results we present, by 
combining reshapings of higher-order tensors into tensors of order three with Proposition 1.32 
in \cite{DDL2013}, as we were kindly informed by its authors. We will nevertheless present an 
example of an identifiable tensor in $S^3 \C^{7}$ of rank $10$ whose identifiability cannot be 
proved by the state-of-the-art specific identifiability criteria from the literature.

It is important to stress that we discuss the general setting of degree $d\ge3$ Veronese embeddings.
We will restrict our attention to nondefective $r$-secants of $\Var{V}$, 
because identifiability will not hold for general tensors on a defective 
$r$-secant variety. This is the interesting setting, 
because the Alexander--Hirschowitz 
theorem \cite{AH1995} stipulates that most $\Sec{r}{\Var{V}}$ are nondefective. 

\subsection{The Hessian criterion}
We recall the main proposition from \cite{COV2014} and adapt it to the present context 
of symmetric tensors.

\begin{lemma0}[Sufficient condition for specific identifiability]\label{infinitepoints}
Let $\Var{V} = v_d(\Pj^n)$ be a nondefective Veronese variety, and let 
$r \le \lceil r_{d,n} \rceil -1$ with $r_{d,n}$ as in \refeqn{eqn_exp_rank}. 
Assume that we are given a \emph{nonsingular} point
\[
 p = p_1 + p_2 + \cdots + p_r \in \Sec{r}{\Var{V}}.
\]
If the linear span of the tangent spaces to $\Var{V}$ at the $p_i$'s, i.e.,
\begin{align*}
 \Plane{M} = \langle \Tang{p_1}{\Var{V}}, \ldots, \Tang{p_r}{\Var{V}} \rangle,
\end{align*}
has the expected dimension, i.e., $r(n+1)$, and if, in addition, the \emph{$r$-tangential contact locus} 
\begin{align*}
 \Var{C}_r = \left\{ p \in \Var{V} \;|\; \Tang{p}{\Var{V}} \subset \Plane{M} \right\}
\subset \Var{V}
\end{align*}
is zero-dimensional at every $p_1, p_2, \ldots, p_r$, then $p$ is
$r$-identifiable, $r$ is its symmetric rank,
and $p = \sum_{i=1}^r p_i$ is its unique decomposition.
\end{lemma0}
\begin{proof}
The proof is obtained by repeating the proofs of \cite[Lemma 4.3, Lemma 4.4, and Theorem 4.5]{COV2014}, 
therein substituting the Segre variety with the Veronese variety $\Var{V}$. 
We present a simplification of the proof of \cite[Theorem 4.5]{COV2014}. 
There is an open neighborhood of $p = p_1 + \cdots + p_r$ consisting of points for which smoothness, Terracini's lemma \cite{Terracini1911}, and the absence of a contact locus will hold. 
The variety must thus be generally identifiable, hence the projection $\pi$ onto 
the first factor of the usual abstract secant
variety $A\Sec{r}{\Var{V}}$ is a birational morphism. After we find another decomposition
 $p =  \sum_{i=1}^r b_iq_i$, we get 
that the fiber $\pi^{-1}(p)$ contains the two points $\left(p, (p_1,\ldots, p_r)\right)$ and 
$\left(p, (q_1,\ldots, q_r)\right)$. Terracini's lemma implies that the connected component of
the fiber passing through $\left(p, (p_1,\ldots, p_r)\right)$ cannot be positive dimensional,
hence, it contains just this unique point.
Since $\pi$ is birational and $p$ is a smooth point of $\Sec{r}{\Var{V}}=
\pi(A\Sec{r}{\Var{V}})$, we have a contradiction with Zariski's Main Theorem.
\end{proof}

\begin{rem0}
 We note a minor omission in the formulation of Theorem 4.5 in \cite{COV2014}, where we forgot 
 to include the condition that $\Plane{M}$ should be of the expected dimension. It is clear 
 from the proof of aforementioned theorem that this condition must hold, as can be understood
 from the invocation of \cite[Lemma 4.3]{COV2014}.
\end{rem0}

\begin{rem0}
 If one chooses $r$ \emph{random} points $p_i \in v_d(\Pj^n)$, then \reflem{infinitepoints} may be 
invoked to prove generic $r$-identifiability. In this way, one can handle the cases
$v_3(\Pj^2)$, $v_3(\Pj^3)$ and $v_3(\Pj^4)$, which were not covered by the proof in 
\refsec{sec:cubics}. The case $v_3(\Pj^1)$ is trivial, because there is only one point, which is 
naturally identifiable.
\end{rem0}

For practically verifying \reflem{infinitepoints}, we need a sufficiently explicit description 
of the $r$-contact locus. This is obtained as follows. Interpreting a point $p_i \in \Var{V}$ 
as a power of a linear form, say 
\[
p_i = ( a_{0,i} x_0 + a_{1,i} x_1 + \cdots + a_{n,i} x_n )^d, 
\]
where $\{x_i\}_{i=0}^n$ is a basis of $\Pj^n$, it follows immediately that the tangent space 
is given by
\[
\Tang{p_i}{\Var{V}} = \langle x_0 ( a_{0,i} x_0 + a_{1,i} x_1 + \cdots + 
a_{n,i} x_n )^{d-1}, \ldots, x_n ( a_{0,i} x_0 + a_{1,i} x_1 + \cdots + a_{n,i} x_n )^{d-1} \rangle. 
\]
If we choose the standard monomial basis 
$\{x_{i_1} x_{i_2} \cdots x_{i_d} \}_{0 \le i_1 \le i_2 \le \cdots \le i_d\le n}$ for $v_d(\Pj^{n})$, 
then this tangent space can be represented in a straightforward manner as a 
$\binom{n+d}{d} \times (n+1)$ matrix of constants, say $T_i$. The Cartesian equations of $\Plane{M}$ 
may then be constructed by computing the kernel of the matrix 
$T = \left[\begin{smallmatrix}T_1 & T_2 & \cdots & T_r \end{smallmatrix}\right]^T$. The number of 
such equations should be precisely $\ell = \binom{n+d}{d} - r(n+1)$; otherwise, the first condition 
in \reflem{infinitepoints} concerning the dimension of the tangent space would be violated. Let us 
denote the Cartesian equations as 
\begin{align}\label{eqn_eqns_in_deal}
 q_l(\vect{x}) = \sum_{i_1=0}^n \sum_{i_2=i_1}^n \cdots \sum_{i_d=i_{d-1}}^n 
k_{(i_1,i_2,\ldots,i_d),l} \cdot x_{i_1} x_{i_2} \cdots x_{i_d} = 0, \quad l = 1, 2, \ldots, \ell,
\end{align}
where the vector $\vect{k}_l = [ k_{(i_1,i_2,\ldots,i_d),l} ]_{0\le i_1 \le i_2 \le \cdots \le i_d\le n}$ 
is the $l$th basis vector of the kernel of the matrix of constants $T$. For imposing that a point 
$\rho = (a_0 x_0 + \cdots + a_n x_n)^d$ is contained in $\Plane{M}$, it should obey the Cartesian 
equations, i.e., $q_l(a_0,a_1,\ldots,a_n) = 0$ for all $l = 1, 2, \ldots, \ell$. That is, the equations \refeqn{eqn_eqns_in_deal} define the ideal-theoretic equations for $\Var{V} \cap \Plane{M}$. 
It similarly follows that deriving the equations \refeqn{eqn_eqns_in_deal} with respect to 
$x_0, x_1, \ldots, x_n$ and substituting $x_0, x_1, \ldots, x_n$ for, respectively, $a_0, a_1, \ldots, a_n$ 
results in the ideal-theoretic equations of the intersection $\Var{C}_r = \Tang{p}{\Var{V}} \cap \Plane{M}$; 
naturally, the $a_0, a_1, \ldots, a_n$ should be treated as new variables. The number of equations thus 
constructed equals $\ell (n+1)$. To determine that $\Var{C}_r$ is zero-dimensional at each $p_i$,
 it suffices to verify that the codimension of the tangent space, i.e., the derivative of the equations 
of the ideal, is $n$ at each of the $p_i$'s. This tangent space can be represented by a matrix $H$ of 
size $(n+1) \times \ell(n+1)$, which contains only constants when it is evaluated at one of the $p_i$'s. The
 rank of $H$ coincides with the dimension of the contact locus and can be computed using simple linear 
algebra. As was remarked in \cite{COV2014}, $H$ can be interpreted as a ``stacked Hessian'' matrix 
$H = \left[\begin{smallmatrix} H^1 & H^2 & \ldots & H^\ell \end{smallmatrix}\right]$, wherein $H^k$ is 
the Hessian matrix of partial derivatives 
\begin{align*}
H^k = [h_{i,j}^k]_{i,j=0}^n = 
\left[ \frac{\partial^2}{\partial x_j \partial x_i} q_k(x_0,x_1,\ldots,x_n) \right]_{i,j=0}^n;
\end{align*}
this is the reason why we call the above approach of verifying \reflem{infinitepoints} the 
\emph{Hessian criterion}.

A computer implementation of the Hessian criterion in Macaulay2 is included in the 
\texttt{specific-identifiability.m2} file that accompanies the arXiv version of this article.

\subsection{The smoothness criterion}
The Hessian criterion in \reflem{infinitepoints} may only 
be applied to smooth points of $\Sec{r}{\Var{V}}$. 
One approach for proving smoothness consists of verifying that the local 
equations of $\Sec{r}{\Var{V}}$ are of the expected degree. Such equations 
are known in the case when the number of terms $r$ in the symmetric 
decomposition is sufficiently small. A standard nontrivial set of 
local equations is generated by the $(r+1)$-minors of the usual symmetric 
flattenings; see \cite[Theorem 7.3.3.3]{Land_book} and 
\cite[Theorems 4.5A and 4.10A]{IK1999}. For Veronese embeddings of odd degree, 
the Young flattenings from \cite[Section 4]{LandsbergOttaviani2013} apply in 
a wider range than the standard symmetric flattenings; however, they are 
more involved to explain and implement. Our discussion will 
focus on the simple symmetric flattenings, which can still handle a respectable 
number of cases for Veronese embeddings of degree at least four. 
For degree three Veronese embeddings, the Young flattenings that were 
described in \cite{LandsbergOttaviani2013} should be employed.

The strategy that we present for proving that $p$ corresponds to a smooth 
point is, essentially, based on \cite[Theorems 4.10A and 4.5A]{IK1999} and 
\cite[Theorem 7.3.3.3]{Land_book}
and consists 
of obtaining local equations of the $r$-secant variety $\Sec{r}{\Var{V}}$ in $p$.  
Crucial to this approach are the \emph{symmetric flattenings}, which, we recall, may be 
defined as follows. Let $p \in S^d \C^{n+1}$, then we can define the map
\begin{align*}
\phi^p_{k}:\qquad  S^k (\C^{n+1})^* &\to S^{d-k} \C^{n+1} \\
 x_{i_1} x_{i_2} \cdots x_{i_k} &\mapsto 
\frac{\partial^k}{\partial x_{i_1} \partial x_{i_2} \cdots \partial x_{i_k}}.
\end{align*}
We have the following.

\begin{lemma0}[Sufficient condition for smoothness]\label{lem_smoothness}
Let $\Var{V} = v_d(\Pj^n)$ be the Veronese variety, let 
$\delta = \lfloor \frac{d}{2} \rfloor$, and let $r < r_{\delta,n}$. Assume that we
are given a point 
\[
 p = p_1 + p_2 + \cdots + p_d \in \Sec{r}{\Var{V}}.
\]
Let $\Plane{N}$ be the following linear space:
\[
 \Plane{N} = 
\operatorname{ker}( \phi^p_{\delta} ) \circ \operatorname{image}( \phi^p_{\delta} )^\perp 
\subset S^d (\C^{n+1})^*,
\]
i.e., the symmetric product of the kernel and the complement of the image of $\phi^p_{\delta}$. If 
\[
\rank{ \phi_\delta^p } = r, \quad\text{and}\quad 
\dim \Plane{N} = \binom{n+d}{d} - r(n+1),
\]
then $p$ is a smooth point of $\Sec{r}{v_d(\Pj^n)}$.
\end{lemma0}
\begin{proof}
The subspace $\Plane{N}$ is the normal space at $p$ of the locus of $(r+1)$-minors of the catalecticant matrix $\phi_\delta^p$;
see, for example, \cite[Prop. 5.3.3.1]{Land_book}.
If $\Plane{N}$ has the
expected dimension $\binom{n+d}{d} - r(n+1)$, then the locus of $(r+1)$-minors of the catalecticant matrix $\phi_\delta^p$ is smooth at
$p$ and of the expected dimension $r(n+1)-1$. The $r$-secant variety $\Sec{r}{v_d(\Pj^n)}$ is contained in that locus, being of the 
expected dimension $r(n+1)-1$ by the Alexander-Hirschowitz theorem, so it too has to be smooth at $p$.
\end{proof}

A natural question concerning the foregoing lemma concerns the maximum value of $r$ for which it can 
be applied. That is, if we pick a sufficiently general smooth point $p \in \Sec{r}{\Var{V}}$, what is 
the maximum value of $r$ for which \reflem{lem_smoothness} can prove that $p$ is, indeed, smooth? A lower bound 
follows immediately from the work of Iarrobino and Kanev \cite[Theorem 4.10A]{IK1999}:
\begin{prop0}\label{prop_ik_lower_bound}
 Let $\Var{V} = v_d(\Pj^n)$ be the Veronese variety, let $\delta = \lfloor \frac{d}{2} \rfloor$, and 
 let 
 \[
  r \le \binom{n+\delta-1}{\delta-1}.
 \]
 Then, \reflem{lem_smoothness} can be applied to all points of an irreducible component of $\Sec{r}{\Var{V}}$
 minus some Zariski-closed set.
\end{prop0}
\begin{proof}
 The claim follows from \cite{IK1999} and the fact that the conditions on the dimension of $\Plane{M}$ 
 and the rank of $\phi_\delta^p$ are valid on dense open sets in the Zariski topology.
\end{proof}

In \reftab{tab_bounds}, some values of the lower bound in \refprop{prop_ik_lower_bound} 
are tabulated along with a sharp maximum value of $r$ 
for which the equations generated by \reflem{lem_smoothness} generate an irreducible component 
of $\Sec{r}{\Var{V}}$. The values of this alleged sharp upper bound were computed by taking 
random points on this variety and verifying \reflem{lem_smoothness}; as such, they are only true 
with high probability. It is clear from the table that the lower bound in \refprop{prop_ik_lower_bound}
is not sharp.\footnote{Recall that the range of applicability of \cite[Proposition 1.32]{DDL2013} may be wider
by combining it with reshapings.}

\begin{table}
\caption{The maximum value $r$ for which \reflem{lem_smoothness} applies to all points in an 
irreducible component of $\Sec{r}{\Var{V}}$ minus some Zariski-closed set is displayed as the middle 
set of columns ($\clubsuit$) for each degree $d=4,5,6,7,8$ of the Veronese embedding 
$\Var{V} = v_d(\Pj^n)$. The left set of columns ($\spadesuit$) shows the lower bound from \refprop{prop_ik_lower_bound}, for every $d$. The right set of columns ($\square$) 
shows the maximum value of $r$ for which Kruskal's criterion is applicable, for every $d$. 
A $\star$ indicates that the value could not be computed within a reasonable time. 
Values displayed in boldface indicate the widest range for $r$ for a particular combination 
of the degree $d$ and size $n$.
}\label{tab_bounds}
\setlength{\tabcolsep}{4pt}
\begin{tabular}{ccccccccccccccccccccccc} 
\toprule
\multirow{2}{*}{$n$} && \multicolumn{19}{c}{$d$} \\
\cmidrule{3-21}
&& \multicolumn{3}{c}{4} && \multicolumn{3}{c}{5} && \multicolumn{3}{c}{6} && \multicolumn{3}{c}{7} &&  \multicolumn{3}{c}{8}\\
\cmidrule{3-5}  \cmidrule{7-9} \cmidrule{11-13} \cmidrule{15-17} \cmidrule{19-21}
&&$\spadesuit$&$\clubsuit$&$\square$&&$\spadesuit$&$\clubsuit$&$\square$&&$\spadesuit$&$\clubsuit$&$\square$&&$\spadesuit$&$\clubsuit$&$\square$&&$\spadesuit$&$\clubsuit$&$\square$ \\
\midrule
1 &&  2 & \bf2  &\bf2 && 2  &    2  &\bf3 && 3   & \bf3  &\bf3 && 3   &    3  &\bf4 && 4    & \bf4  &\bf4\\
2 &&  3 & \bf4  &\bf4 && 3  &    4  &\bf5 && 6   & \bf6  &\bf6 && 6   & \bf7  &\bf7 && 10   & \bf10 & 8\\
3 &&  4 &    5  &\bf6 && 4  &    6  &\bf8 && 10  & \bf12 & 9   && 10  & \bf15 & 11  && 20   & \bf23 & 12\\
4 &&  5 &    7  &\bf8 && 5  &    9  &\bf10&& 15  & \bf21 & 12  && 15  & \bf27 & 14  && 35   & \bf47 & 16\\
5 &&  6 & \bf10 &\bf10&& 6  & \bf14 & 13  && 21  & \bf33 & 15  &&\bf21&$\star$& 18  && 56   & \bf87 & 20\\
6 &&  7 & \bf12 &\bf12&& 7  & \bf19 & 15  && 28  & \bf50 & 18  &&\bf28&$\star$& 21  &&\bf84 &$\star$& 24\\
7 &&  8 & \bf16 & 14  && 8  & \bf25 & 18  && 36  & \bf72 & 21  &&\bf36&$\star$& 25  &&\bf120&$\star$& 28\\
8 &&  9 & \bf20 & 16  && 9  & \bf33 & 20  &&\bf45&$\star$& 24  &&\bf45&$\star$& 28  &&\bf165&$\star$& 32\\
9 && 10 & \bf25 & 18  && 10 & \bf41 & 23  &&\bf55&$\star$& 27  &&\bf55&$\star$& 32  &&\bf220&$\star$& 36\\
10&& 11 & \bf29 & 20  && 11 &$\star$&\bf25&&\bf66&$\star$& 30  &&\bf66&$\star$& 35  &&\bf286&$\star$& 40\\
\bottomrule
\end{tabular}
\end{table}

For $d=3$, the symmetric flattenings are only sufficient for $r=1$ and $2$. One should employ 
Young flattenings \cite{LandsbergOttaviani2013} for extending the range instead. As an illustration of the range that can be covered by such equations, we present in \reftab{tab_third_order} the maximal value of $r$ for which the equations generated by Young flattenings generate an irreducible component of $\Sec{r}{v_3(\Pj^n)}$. The range of applicability is also compared with the criterion of Domanov and De Lathauwer \cite[Proposition 1.32]{DDL2013}, which to the best of our knowledge provides the state-of-the-art specific identifiability criterion.
The border case $\Sec{11}{v_3(\Pj^6)}$ in \reftab{tab_third_order}
is special, because in all of the random choices of $p\in \Sec{11}{v_3(\Pj^6)}$ that we tested, the Young flattening had rank two less than expected. Still, the corresponding minors of the Young flattening did cut $\Sec{11}{v_3(\Pj^6)}$ scheme-theoretically at $p$.

\begin{table}
\caption{The maximum value $r$ for which Young flattenings generate an
irreducible component of $\Sec{r}{\Var{V}}$ is given in the first row ($\clubsuit$). 
The second row ($\diamondsuit$) shows the maximum value of $r$ for which the state-of-the-art 
criterion of \cite[Proposition 1.32]{DDL2013} is applicable. 
Values displayed in boldface indicate the widest range for $r$ for a particular combination 
of the degree $d$ and size $n$.} \label{tab_third_order}
\begin{tabular}{cccccccccc}
\toprule
$n$            & 1    & 2    & 3    & 4    & 5    & 6     & 7     & 8     & 9     \\
\midrule
$\clubsuit$    & \bf2 & \bf3 & \bf5 & \bf6 & \bf8 & \bf11 & \bf11 & \bf14 & \bf15  \\
$\diamondsuit$ & \bf2 & \bf3 & \bf5 & \bf6 & \bf8 & 9     & \bf11 & 12    & 14     \\
\bottomrule
\end{tabular}
\end{table}

An implementation in Macaulay2 of the above sufficient condition for smoothness based on a 
symmetric flattening is included in the \texttt{specific-identifiability.m2} file that is 
provided with the arXiv version of this paper. In addition, this file contains an
 implementation of a smoothness test based on Young flattenings for Veronese embeddings of 
degree three.

\subsection{An elementary algorithm}\label{sec_algorithm}
For the sake of completeness, we present an algorithm that attempts
to prove the identifiability of a given Waring decomposition 
$p = p_1 + p_2 + \cdots + p_r \in \Sec{r}{\Var{V}}$, where $\Var{V} = v_d(\Pj^n)$, 
by checking the sufficient conditions in \reflem{infinitepoints} and \reflem{lem_smoothness}. 
It operates as follows.
\begin{enumerate}
 \item[S1.] Construct a matrix representation of the span 
$\Plane{M} = \langle \Tang{p_1}{\Var{V}}, \Tang{p_2}{\Var{V}}, \ldots, \Tang{p_r}{\Var{V}} \rangle$. 
If $\rank \Plane{M} < r(n+1)$, then the algorithm terminates, claiming 
that it cannot prove the identifiability of $p$.
 
 \item[S2.] Construct the symmetric flattening $\phi^p_\delta$ for $\delta = \lfloor \frac{d}{2} \rfloor$. 
If $\operatorname{rank} \phi^p_\delta < r$, then the algorithm terminates, claiming that it 
cannot prove identifiability of $p$.

 \item[S3.] Compute the matrix $\Plane{N} = \operatorname{ker}( \phi^p_{\delta} ) \circ 
\operatorname{image}( \phi^p_{\delta} )^\perp$. If $\rank \Plane{N} > \binom{n+d}{d} - r(n+1)$, 
then the algorithm terminates, claiming that it cannot prove identifiability of $p$.

 \item[S4.] Compute a basis of the kernel of $\Plane{M}$. Denote the number of equations by $\ell$.

 \item[S5.] For every point $p_i$, $i = 1, 2, \ldots, r$, perform the following:
\begin{enumerate}
 \item[S5a.] Construct the Hessians $H^k$ for $k = 1, 2, \ldots, \ell$ evaluated at the point $p_i$ 
and stack them into the matrix $H$.
 \item[S5b.] If $\operatorname{rank} H < n$, then the contact locus is of positive dimension at $p_i$. 
The algorithm halts, claiming that it cannot prove identifiability of $p$.
\end{enumerate}

 \item[S6.] The algorithm proclaims that the Waring decomposition $p = p_1 + \cdots + p_r$ is unique 
and that $p$ is a smooth point of $\Sec{r}{\Var{V}}$.
\end{enumerate}

It is instructive to investigate the largest value of $r$ for which the above algorithm may be 
expected to prove identifiability of a sufficiently general point $p \in \Sec{r}{\Var{V}}$ on
a generically identifiable Veronese variety $\Var{V}$. As the Hessian criterion applies for all 
tensors of subgeneric rank, it follows that the range of applicability is bounded only by 
the smoothness test. That is, the highlighted columns in \reftab{tab_bounds} contain the relevant
values. A popular criterion for testing identifiability that 
is applicable for Veronese embeddings of any degree $d\ge3$ is the so-called Kruskal condition 
\cite{Kruskal1977,JS2004}. Let $p_i = \vect{a}_i^{\otimes d} \in \Var{V}$ be some specific points.
In the symmetric setting, Kruskal's condition states that if
\[
 r \le \frac{1}{2}( d k - d + 1 ),
\]
where $k$ is the largest number such that every subset of $\{\vect{a}_i\}_i$ consisting of $k$ vectors 
is linearly independent. For points in general configuration, the maximum value for $k$ is thus $n+1$.
A comparison between the proposed criterion for specific identifiability and Kruskal's criterion is 
also featured in \reftab{tab_bounds}. 

\subsection{The algorithm at work for a specific example}
Consider the following $10$-term Waring decomposition in $S^3 \C^7$:
\begin{align*}
 p &= \sum_{i=0}^6 x_i^3 + \\
&\quad (4x_0 + 3x_1 + 2x_2 + x_3 + 2x_4 + 3x_5 + 4x_6)^3 + \\
&\quad (x_0 + x_1 + 2x_2 + 2x_3 + 3x_4 + 3x_5 + 4x_6)^3 \\
&\quad (x_0 + 2x_1 + 3x_2 + 4x_3 + 5x_4 + 6x_5 + 7x_6)^3.
\end{align*}
The identifiability of this example cannot be handled with Domanov and De Lathauwer's criterion,
because it only applies for Waring decompositions with $r \le 9$. The 
sufficient condition presented in this paper, on the other hand, is applicable up to $11$ 
terms. Therefore, we can run the algorithm presented above. This example may be verified with 
the \texttt{specific-identifiability.m2} script that is provided with the arXiv version 
of this paper. In the first step, the $84 \times 70$ matrix representing the span of 
$\Plane{M}$ is constructed. Its rank is $70$, as expected. The algorithm proceeds with the 
construction of Young flattening, checking they generate an irreducible component. Then,
a basis of the kernel of $\Plane{M}$ is computed, containing $14$ equations. Note that the 
dimension of $\Plane{N}$ and the codimension of $\Plane{M}$ should always be equal in the 
approach for certifying identifiability that was proposed in this section. Then, for each of
the points, the $7\times7$ Hessian matrices are computed and stacked into a $7 \times 98$ 
matrix $H$. The rank of $H$ equals $6$ for each of the points. In addition, the kernel of
$H^T$ consists of a single vector that must be a multiple of the coefficient vector of $p_i$; 
for example, the vector in the kernel of the stacked Hessian $H$ corresponding to the last 
term in the Waring decomposition of $p$ is a multiple of its coefficient vector 
$\left[\begin{smallmatrix}1 & 2 & 3 & 4 & 5 & 6 & 7 \end{smallmatrix}\right]$. 
For each of the points, this is indeed the case. Finally, the algorithm positively concludes
that $p$ admits a unique Waring decomposition, i.e., $p$ is identifiable.

\end{document}